\documentclass[a4paper]{amsart}
\usepackage{amssymb,amsmath,amsthm,amscd,graphicx}

\def \op{\operatorname}
\def \lra{\longrightarrow}

\def \wh{\widehat}

\newtheorem{theorem}{Theorem}
\theoremstyle{plain}

\newtheorem{proposition}[theorem]{Proposition}

\newcommand{\thmref}[1]{Theorem~\ref{#1}}

\newcommand{\propref}[1]{Proposition~\ref{#1}}

\title [Three manifolds as geometric branched coverings]{Three manifolds as geometric branched coverings of the three
sphere}
\author [Brumfiel,  Hilden,  Lozano,  Montesinos,
Ramirez, Short, Tejada,  Toro]{G.~Brumfiel, H.~Hilden,
M.T.~Lozano*, J.M.~Montesinos--Amilibia*, E.~Ramirez--Losada,
H.~Short, D.~Tejada, M.~Toro}
\address[G.~Brumfiel]{Stanford University, Stanford, Ca, USA}
\email[G.~Brumfiel]{brumfield@math.stanford.edu}
\address[H.~Hilden]{University of Hawaii, Honolulu, Hi, USA}
\email[H.~Hilden]{mike@math.hawaii.edu}
\address[M.T.~Lozano]{Universidad de Zaragoza, Zaragoza, Spain}
\email[M.T.~Lozano]{tlozano@unizar.es}
\address[J.M.~Montesinos]{Universidad Complutense, Madrid, Spain}
\email[J.M.~Montesinos]{montesin@mat.ucm.es}
\address[E.~Ramirez]{CIMAT, Mexico} \email{kikis@cimat.mx}
\address[H.~Short]{Universite de Provence, Marseille, France}
 \email{hamish.short@cmi.uni-mrs.fr}
 \address[D.~Tejada, M.~Toro]{Universidad Nacional de Colombia, Madellin, Colombia}
 \email{dtejada@unalmed.edu.co, mmtoro@unalmed.edu.co}
\thanks {$^*$This research
was supported by grants MTM2004088080 and MTM2006-00825}
\begin{document}
\begin{abstract}
One method for obtaining every closed orientable 3-manifold is as
a branched covering of $S^3$ branched over a link.   There are
topological results along these lines that cannot be improved upon
in two respects:
 (\emph{A}), The minimum possible number of sheets in the covering is
 three;
(\emph{B}), There are individual knots and links (universal knots
and links) that can serve as branch set for every 3-manifold,
$M^3$. Given the growing importance of geometry in 3-manifold
theory it is of interest to obtain geometrical versions of
topological results (\emph{A}) and (\emph{B}). Twenty years ago a
geometric version of result (\emph{B}) was obtained using
universal groups. In this paper we obtain a geometric version of
result (\emph{A}), also by means of universal groups.
\end{abstract}

\maketitle

\begin{section}{Introduction}

Some time ago two of the co--authors of this paper showed that all
closed orientable 3--manifolds are threefold irregular simple
branched coverings of $S^3$ with branch set a knot of a link.
(\cite{Hilden1974},\cite{Hirsch1974},\cite{Montesinos1974}).  Both
the statement and proof of this theorem were purely topological in
character. Geometry played no role.

The central purpose of this paper is to introduce the idea of a
geometric branched covering and then to show that every closed
orientable 3--manifold $M^3$ is a threefold irregular simple
branched covering of $S^3$ that is ``geometric" in the following
sense:

There is a universal group $U$; the orbifold group of the Borromean
rings with singular angle ninety degrees.  Thus $U$ is a finite
covolume Kleinian group of hyperbolic isometries of $H^3$.

There are finite index subgroups of $U$, $G$, and $G_1$, such that
$M^3 = H^3/G$, $S^3 = H^3/G_1$, and $G$ is a non normal index three
subgroup of $G_1$.

The groups $U$, $G$, and $G_1$ induce hyperbolic orbifold
structures on $S^3$, $M^3$, and $S^3$ respectively and the maps
induced by group inclusions, \[ M^3 = H^3/G \lra S^3 = H^3/G_1
\lra S^3 = H^3/U \, ,\] are branched covering space maps with
$M^3\lra S^3$ being three to one.

The new theorem, \thmref{t9} of this paper, can be considered to be
the ``geometrization" of the old theorem.

This work is closely related to the concepts of universal group and
universal link.

A knot or link is said to be universal if every closed orientable
$3$--manifold occurs as a finite branched covering space of $S^3$
with branch set the knot or link. W.~Thurston introduce this
concept in his paper \cite{Thurston1982}, where he also exhibited
some universal links, and asked if the Whitehead link and the
Figure eight knot were universal. In \cite{HLM1983c},
\cite{HLM1985a} three of the co-authors answered Thurston's
question in the affirmative and proved that every non toroidal
rational knot or link is universal, and that the Borromean rings
are universal. Subsequently many other universal knots and links
were found: \cite{Uchida1991}, \cite{Uchida1992},
\cite{Nunez2004}, \cite{HLM1985b}, \cite{HLM1987}, \cite{LM1997},
\cite{HLM2004}.

Given an n-tuple of classical knots in $S^3$
$(K_{1},\,K_{2},...K_{n})$ such that $K_{i}$ doesn't intersect
$K_{j}$ for $i\neq j$  and an $n$-tuple of positive integers
$(m_{1},m_{2},...m_{n})$ we say that the link $L$ that is the union
of the knots is $(m_{1},m_{2},...m_{n})$-\emph{universal} if every
closed oriented 3-manifold $M$ is a branch cover of $S^3$ with
branch set $L$ such that for all $j$ the branch index at a component
of the preimage of $K_{j}$ is an integer, possibly one, dividing
$m_{j}$. A $(n,n,...n)$-universal link will be called $n$-universal
for short. A $(m_{1},m_{2},...m_{n})$-universal link $L$ is \emph{
minimal} if $m_{i}>1$, for all $i$, and there is no $n$-tuple
$(j_{1},j_{2},...j_{n})$, not equal $(m_{1},m_{2},...m_{n})$, with
$j_{i}$ dividing $m_{i}$ for all $i$ for which $L$ is
$(j_{1},j_{2},...j_{n})$-universal.

Closely related to the concept of universal knot or link is that of
universal group.  A finite covolume, discrete group of hyperbolic
isometries $U$, acting on $H^3$, is said to be universal if every
closed orientable $3$--manifold $M^3$ occurs as a quotient space,
$M^3 = H^3/G$, where $G$ is a finite index subgroup of $U$.  Such
groups $U$ must contain rotations, else all $3$--manifolds including
$S^3$ would have hyperbolic structure.

Also, given $M^3$, there are infinitely many finite index $G$'s with
$M^3 = H^3/G$ so that this doesn't give anything like a
classification of $3$--manifolds.  It can be considered analogous to
Heegaard splittings or Kirby calculus presentations.  We know that
$S^3 = H^3/U$ so that $S^3$ inherits a hyperbolic orbifold structure
from $U$.

Some of the co--authors were involved in the proof of the existence
of a universal group $U$ (\cite{HLMW1987}).  The group defined in
(\cite{HLMW1987}), which from now on we denote by $U$, is the
orbifold group of the Borromean rings with ninety degrees.  This
group $U$ thus induces a hyperbolic orbifold structure on $S^3$ with
singular set the Borromean rings, and singular angle ninety degrees.
Unfortunately the proof of the universality of $U$ in
(\cite{HLMW1987}) cannot be adapted to prove the geometric branched
covering space theorem referred to earlier.

There is a new proof that serves our purposes.  It starts out
following (\cite{HLMW1987}), then follows ideas in
(\cite{HLM1983c}), then follows section 5 of (\cite{LM1997}) in
which infinitely many 2-universal links are defined, and then uses a
branched covering of $S^3$ by $S^3$ with branch set the Borromean
rings and branching indexes \{1,2\} ( branching of type \{1,2\} for
short) such that the ``doubled" Borromean rings (2-universal) occur
as a sublink of the preimage of the Borromean rings.

Rather than put the reader through the difficult task of actually
finding these references in some library we prefer to give a new,
relatively self--contained proof of the universality of $U$.

In the next section we state and prove the geometric branched
covering space theorem.
\end{section}

\begin{section}{ Geometric branched covering space theorem}

Our point of departure is the following theorem
(\cite{Hilden1974},\cite{Hirsch1974},\cite{Montesinos1974}).

\begin{theorem}\label{m31} Let $M^3$ be a closed orientable $3$--manifold. Then $M^3$ is
a 3 to 1 irregular simple branched cover of $S^3$ with branch set a
knot or link (as opposed to a graph).
\end{theorem}

Such 3 to 1 coverings with branch set a link $L$ correspond to
transitive representations $\rho: \pi_1(S^3-L)\lra \Sigma_3$ in
which meridians are sent to transpositions.  Years ago Ralph Fox
(\cite{CrowellFox}) had the genial idea of representing
transpositions by colours.  Throughout the paper we shall follow
this idea; $\op{Red}=R=(12)$, $\op{Yellow}=Y=(23)$,
$\op{Blue}=B=(13)$.

Then, given a classical knot or link diagram, simple transitive
representations to $\Sigma_3$ (Simple means meridians go to
transpositions), correspond to colourings of the bridges such that
at each crossing either all three bridges are the same colour or all
three have different colours and at least two colours are used.
(This itself is equivalent to the Wirtinger relations, which have
form $xy=yz$, being satisfied.)

Given a 3--1 simple branched covering $p: M^3 \lra S^3$ branched
over the coloured link $L$, there is a move, illustrated in Figure 1
and called a Montesinos move \cite{Montetesis}, which changes the
coloured link $L$ to a different coloured link $L'$.  The change
takes place inside a ball.  Although $L$ is changed, the topological
type of $M^3$ is not.  The reason is that the $3$--fold simpler
cover of a ball branched over two unknotted, unlinked arcs is a
ball.  Thus doing a Montesinos move on a link in $S^3$ is equivalent
to removing a ball from $M^3$ and sewing it in differently.  We call
a sequence of Montesinos moves a Montesinos transformation.

\begin{figure}
[ptbh]
\begin{center}
\includegraphics
{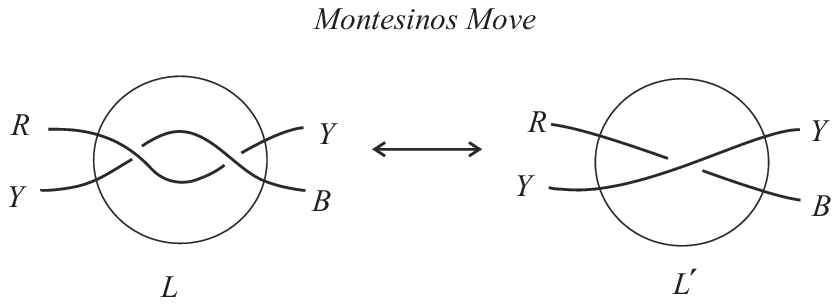}%
\caption{}%
\end{center}
\end{figure}

We give several examples of transformations of links using
Montesinos moves and isotopies which will be useful to us in the new
proof of the universality of $U$.

\begin{figure}
[ptbh]
\begin{center}
\includegraphics
{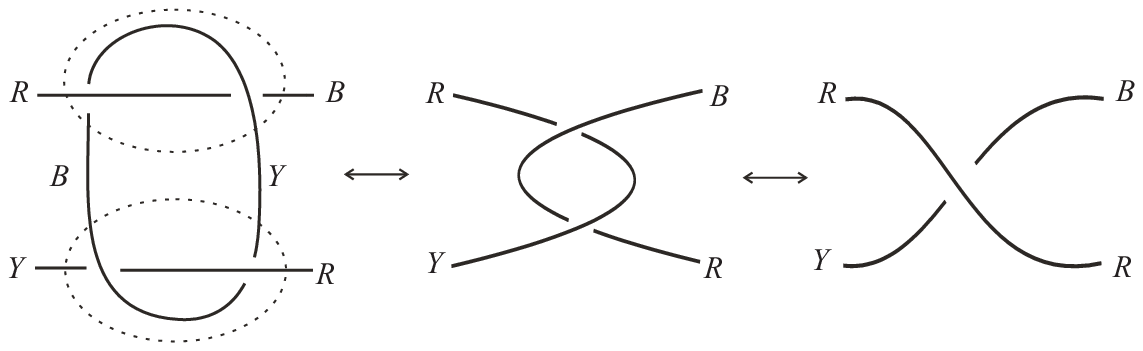}%
\caption{}%
\end{center}
\end{figure}

\begin{figure}
[ptbh]
\begin{center}
\includegraphics
{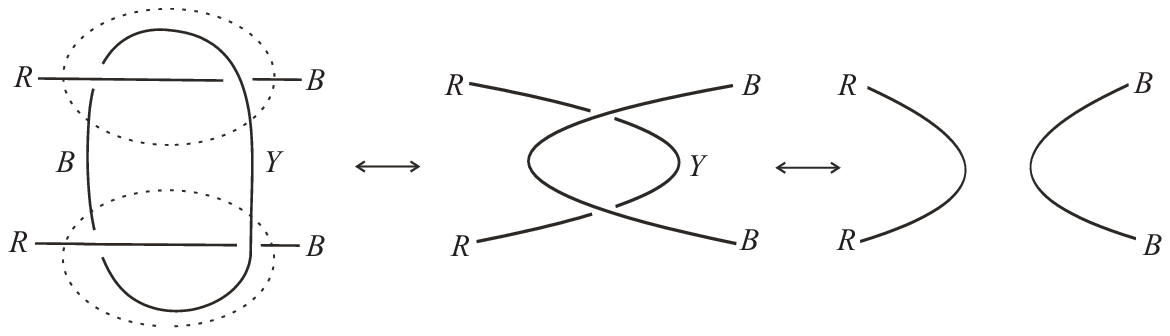}%
\caption{}%
\end{center}
\end{figure}

\begin{figure}
[ptbh]
\begin{center}
\includegraphics
{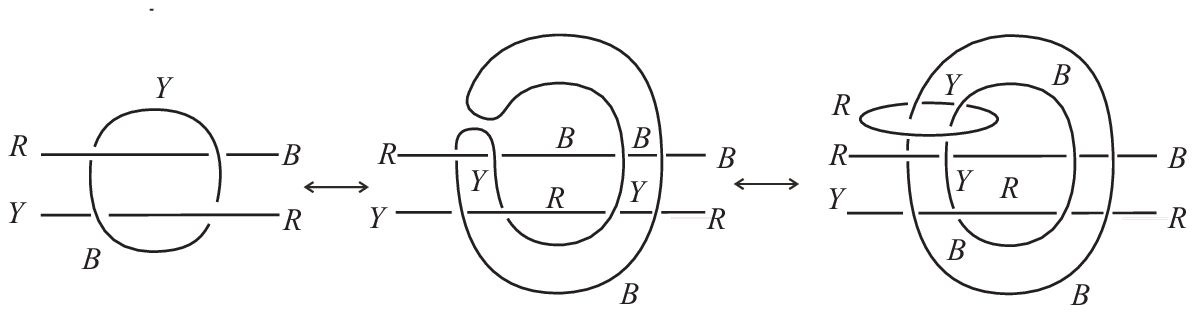}%
\caption{}%
\end{center}
\end{figure}

\begin{figure}
[ptbh]
\begin{center}
\includegraphics
{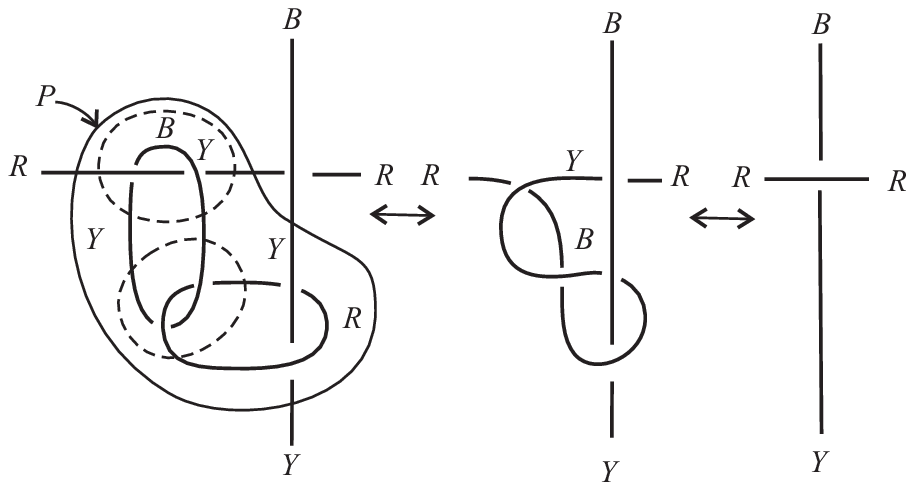}%
\caption{}%
\end{center}
\end{figure}

\begin{figure}
[ptbh]
\begin{center}
\includegraphics
{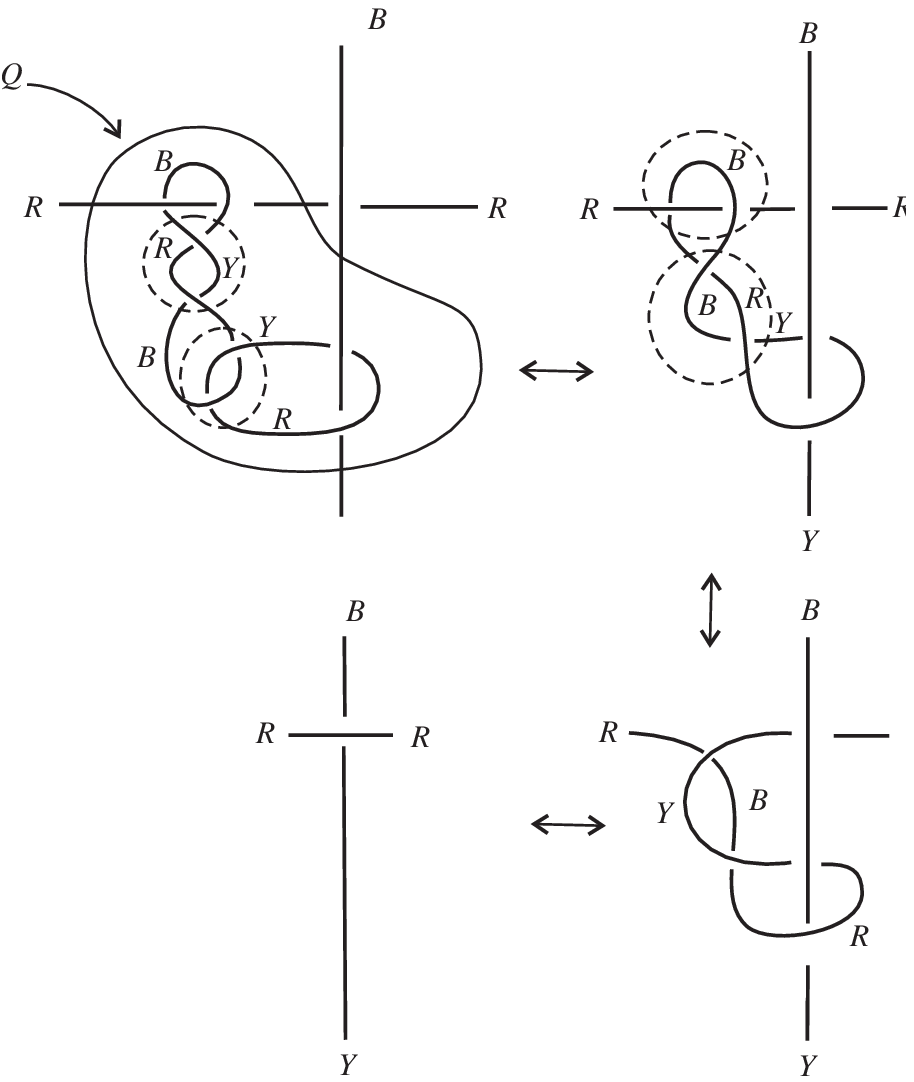}%
\caption{}%
\end{center}
\end{figure}

Now we are ready to prove the universality of $U$.  Let $M^3$ be a
closed orientable \break$3$--manifold and let $p:M^3\lra S^3$ be a
simple $3$--fold covering branched over the link $L$.  We shall
apply a Montesinos transformation to $L$ to obtain a link $L'$ which
suits our purposes.  For definitiveness we shall think of $S^3$ as
$E^3\cup\{\infty\}$, assume $L$ is contained in $E^3$ and use
cylindrical coordinates in $E^3$.

As every link is a closed braid we can assume $L$ is a closed braid.
This means that each component of $L$ has a parametrization $(r(t),
\theta(t), z(t))$ in which $\theta(t)$ is strictly increasing and
the projection on the plane $z=0$ is ``nice".  (i.e., there are no
triple points.)

Using an isotopy of the type illustrated in Figure 7 we can assume
every crossing has 3 colors.

\begin{figure}
[ptbh]
\begin{center}
\includegraphics
{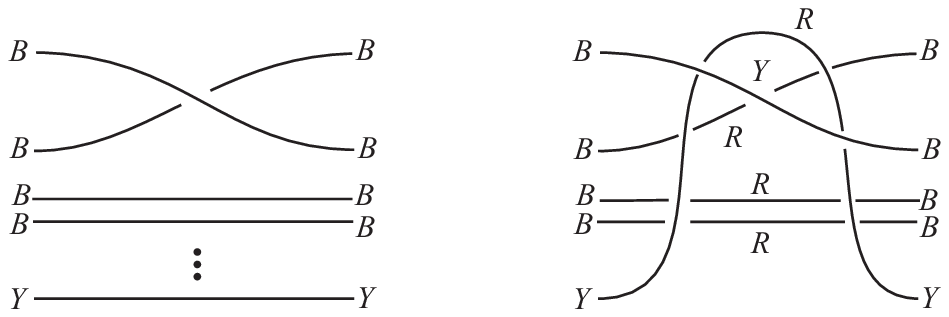}%
\caption{}%
\end{center}
\end{figure}

Then, using Montesinos moves as in Figure 1, we can assume all
crossings are ``positive" (See Figure 8)

\begin{figure}
[ptbh]
\begin{center}
\includegraphics
{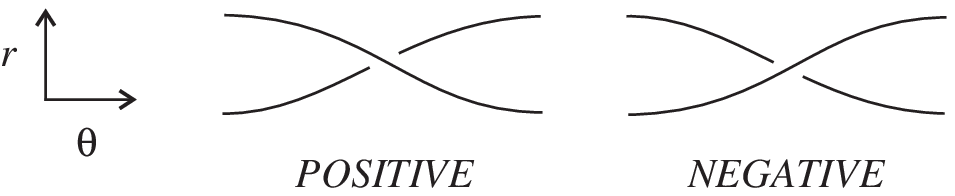}%
\caption{}%
\end{center}
\end{figure}

We replace each crossing with a new small circle component, as in
the left hand side of Figure 2.  After an isotopy we can assume
that our link $L$ has two types of components; ``braid" or
``horizontal" components that lie in the plane $z=0$ and have
equations $r=\op{constant};\theta=\op{arbitrary}, z=0$, and
``small circle" components, whose projection on the plane $z=0$
appears as in the left side of Figure 2.

Next, using the Montesinos transformation of Figure 4, (which is
validated by using the Montesinos transformation of Figure 3.), we
replace each small circle component by three components as in the
right hand side of Figure 4.  Now our link $L$ has three types of
components; horizontal components, big circle components and small
circle components as illustrated in the right hand side of Figure 4.
Each small circle component links two big circle components.

We isotope out link $L$ so that each big circle component appears as
in the left hand side of Figure 9, that is, it extends over the top
and bottom of all the horizontal components.

\begin{figure}
[ptbh]
\begin{center}
\includegraphics
{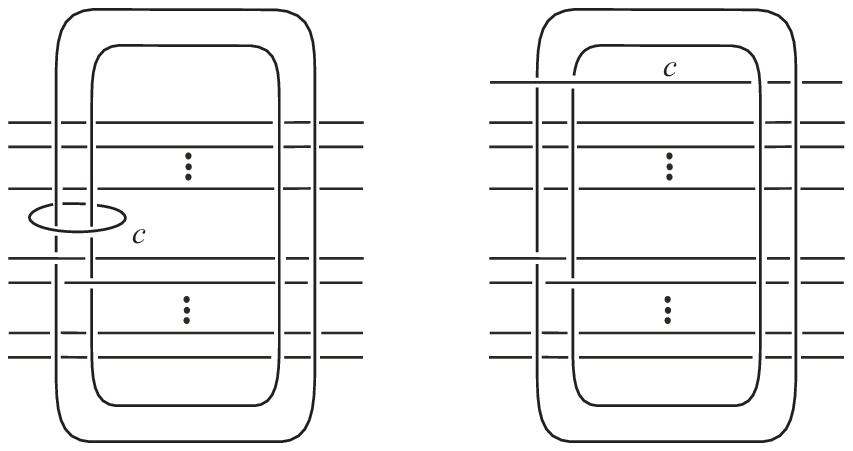}%
\caption{}%
\end{center}
\end{figure}

Then we isotope the small circle components, one at a time, so that
they become braid or horizontal components.  As we do this to a
particular small circle component ``$c$" as in the left hand side of
Figure 9 it becomes, briefly, the topmost braid component.

Now our link $L$ has two types of components, horizontal components
which lie in the plane $z=0$ and have equation of form
$r=\op{constant}$, $\theta=\op{arbitrary}$, and $z=0$ and large
circle components, which we now call vertical components, whose
projections on the plane $z=0$ are rectangles in the $(r,\theta)$
coordinate system.  Crossings are called horizontal or vertical as
indicated in Figure 10.

\begin{figure}
[ptbh]
\begin{center}
\includegraphics
{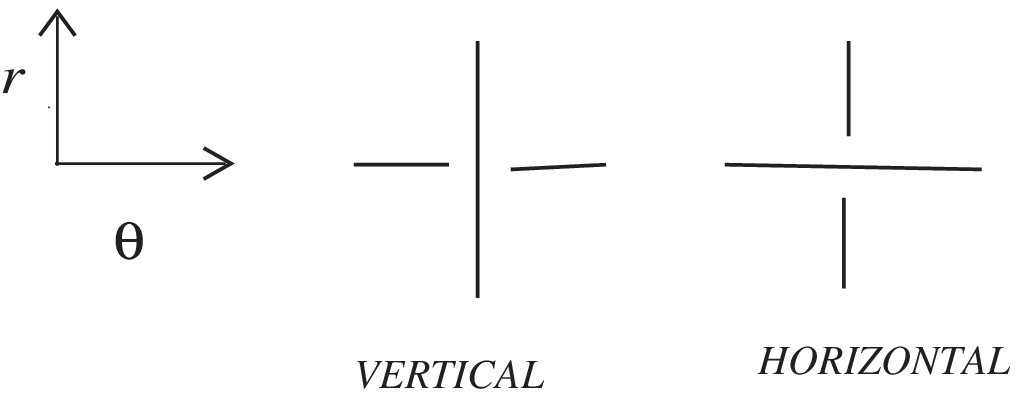}%
\caption{}%
\end{center}
\end{figure}

We observe, and this is crucial in what follows, that every
horizontal crossing is \break$3$--colored.  (The vertical crossings
may not be.)  Also, a particular vertical component links a
particular horizontal component if and only if the left crossing on
the right hand side of Figure 9 is horizontal.  The next step will
be to eliminate horizontal crossings.  We refer the reader to
Figures 5 and 6.

We use the Montesinos transformations illustrated in either Figures
5 or 6 (both are useful) to replace each horizontal crossing by a
vertical crossing.  In the course of doing this, new components, as
indicated in the left sides of Figures 5 and 6, are introduced.
These are contained in the ``peanut shaped" balls indicated by a
``$P$" or ``$Q$" in Figures 5 and 6.

After a slight isotopy our link $L$ has three types of components;
horizontal (with equation $r=\op{constant}$, $z=0$,
$\theta=\op{arbitrary}$); vertical (lying in the plane
$z=\varepsilon >0$, whose projection in the plane $z=0$ is a
rectangle in $(r,\theta)$ coordinates); and ``special".  Each
special component is contained in a ``peanut shaped" topological
ball lying in the region $-\varepsilon \leqq z \leqq 2\varepsilon $
and having projection on the $z=0$ plane as indicated in either the
left or right hand sides of Figure 11.  (We can use one or the other
but never both in the same proof.)

\begin{figure}
[ptbh]
\begin{center}
\includegraphics
{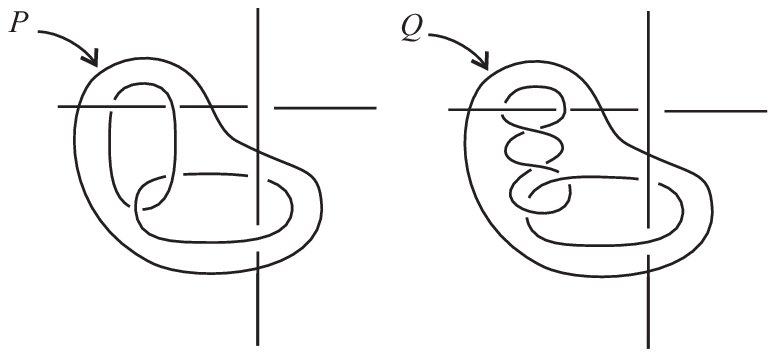}%
\caption{}%
\end{center}
\end{figure}

Thus we have shown that every closed orientable $3$--manifold $M^3$
is a $3$--fold simple branched covering of a link of very special
type.  We shall isotope the link some more and then state a theorem.
We shall use cylindrical coordinates in $S^3=E^3\cup\{\infty\}$ to
better describe the link.  At the moment our link is contained in
the region $[-\varepsilon \leqq z \leqq 2\varepsilon,\, 0 < r_0<r_1
\leqq R_0,\, \theta = \text{arbitrary}]$

We isotope the region $-\varepsilon \leqq z \leqq 2\varepsilon,
r_0 \leqq r \leqq R_0$ without changing the $\theta$ coordinate,
so that the link lies in the thickened toroidal region $1- \delta
\leqq z^2+(r-2)^2\leqq 1+\delta$ and so that the horizontal
components lie in the torus $z^2+(r-2)^2=1$, have equations of
form $[r=\op{constant}, z=\op{constant}, \theta=\op{arbitrary}]$,
and are ``evenly spaced".  (This means that after the isotopy the
images of the horizontal components intersect the circle
$z^2+(r-2)^2=1,\theta=\theta_0$, in $n$ evenly spaced points on
the circle.)

Next we isotope the vertical components so that their images lie
on the torus $z^2+(r-2)^2=1-\delta$, have equations
$\theta=\op{constant}$, $z^2+(r-2)^2=1-\delta$, and are ``evenly
spaced" which means that if there are $m$ vertical components the
constants referred to are $\{2\pi j/m; 0\leqq j\leqq m-1\}$.

It is convenient to introduce new ``toroidal" coordinates $(\rho,
\theta, \varphi)$ which are well defined in a neighborhood of the
torus $z^2+(r-2)^2=1-\delta$.  These are given by equations
$\theta=\theta$; $\rho = \sqrt{z^2+(r-2)^2} =$ distance of a point
from the circle $[z=0 ,\, r=2, \, \theta=\op{arbitrary}]$; $\sin
\varphi = z/\rho$ and $\cos \varphi = r-2/\rho$.  Using these
``toroidal" coordinates we can define projections of $E^3 - \{(z
\op{axis}) \cup (\op{circle}\ [r=2; z=0]\}$ onto the torus
$z^2+(r-2)^2=1$ or the torus $z^2+(r-2)^2=1-\delta$ by $(\rho,
\theta, \varphi) \lra (\theta, \varphi)$.

We summarize the results of the preceding isotopies in a theorem.

\begin{theorem}\label{t1} Let $M^3$ be a closed oriented $3$--manifold. Then
there is a $3$--fold simple branched covering $p:M^3\lra S^3$
branched over a link $L$.

The link $L$ has three types of components.

a. Horizontal.  These lie in the torus $\rho=1$ or $z^2+(r-2)^2=1$
and have equations $[\theta=\op{arbitrary}; \varphi=2 \pi j /n$,
$0\leqq j \leqq n-1, \rho=1]$

b.  Vertical.  These lie in the torus $\rho =.99$ or
$z^2+(r-1)^2=(.99)^2$ and have equations $\varphi=\op{arbitrary},
\theta=2\pi j/ m, 0\leqq j \leqq m-1, \rho=.99$

c.  Special.  These have local projections on the torus $\rho=1$ as
in either the left or right hand side of Figure 11.  The vertical
coordinate is $\varphi$, the horizontal coordinate is $\theta$.
\end{theorem}

In this proof we will use all left hand side or all right hand
side of Figure 11. Both are useful.

Now we find it useful to define two rotations $T_1$ and $T_2$ of
$S^3=E^3\cup \{\infty\}$.  The rotation $T_1$ is simply the
$m$--fold rotation about the $z$--axis given by $T_1:(r,\theta,z)
\lra (r, \theta+2 \pi/m,z)$ the rotation $T_1$ leaves invariant
the set of horizontal and the set of vertical components of the
link $L$.  The rotation $T_2$ is more difficult to describe in
coordinates and we shall not attempt to do so.  Instead we
indicate its important properties.  The rotation $T_2$ has as its
axis the circle $z=0$, $r=2$, $\theta =$ arbitrary; it has order
$n$ and leaves the $\theta$ coordinate unchanged.  It leaves the
set of horizontal and the set of vertical components of $L$
invariant.  It cyclically permutes the horizontal components and
it sends each vertical component to itself.  Its restriction to a
vertical component is just the usual $n$--fold rotation of a
circle.

At this point we must decide whether to use the branch set in the
left or right hand side of Figure 11.  We choose the left for
purposes of illustration.

Using the projection $(\rho, \theta, \varphi)\lra (\theta, \varphi)$
a portion of the image of the link $L$ appears as in Figure 12.
Some of the ``peanut shaped" balls contain two component links and
arcs from a vertical and horizontal component, others contain only
arcs from a vertical and horizontal component.

\begin{figure}
[ptbh]
\begin{center}
\includegraphics
{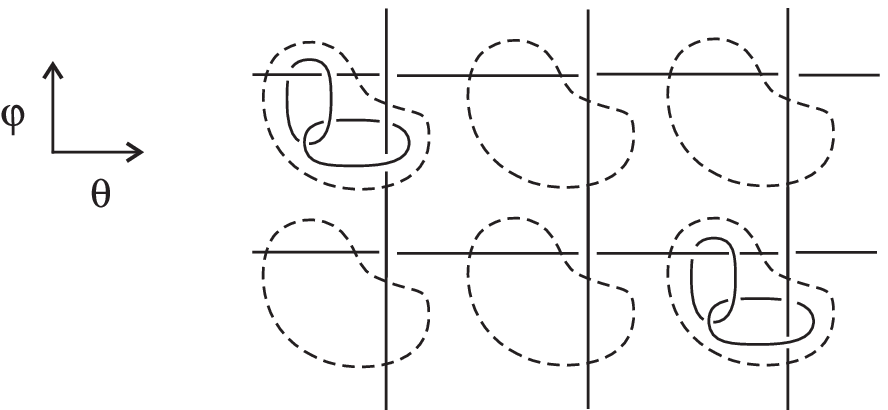}%
\caption{}%
\end{center}
\end{figure}

We can and shall assume that both rotations $T_1$ and $T_2$ leave
the ``expanded" link $L$ invariant.

Now consider the map $f_1: S^3\lra S^3/T_1=S^3$ which is an
$m$--fold cyclic branched covering $S^3\lra S^3$ with branch set the
trivial knot or $z$--axis.  The branch set for the composite map
$f_1\circ p:M^3\lra S^3$ consists of the union of the branch set for
$f_1$ and the image $f_1$ (branch set for $p$) $= z$--axis $\bigcup
f_1(L)$.

The branching is of type $\{1,2\}$ on $f_1(L)$ and of type $\{m\}$
on the $z$--axis$\cup \{\infty\}$.  This means that any disc which
belongs to the preimage of a meridian disk $D$ for $f_1(L)$ is
either mapped homeomorphically $(1-1)$ to $D$ or mapped to $D$ as
a $2-1$ branched covering.  The preimage of a meridian disk $D$
for the $z$--axis is mapped to $D$ as an $m$--fold branched
covering.  The overall branching is of type $\{1,2,m\}$ for the
map $f_1\circ p$ and the $``m"$ is bad for our purposes.  Shortly
we shall show how to change the map $f_1$ to a map $f$ so that the
branching for $f\circ p$ is of type $\{1,2\}$ but first we observe
that the part of the branch set $f_1(L)$ for $f_1\circ p$ contains
$n$ horizontal components, $n$ ``peanut" components but only one
vertical component.

Consider the connected $k$--fold branched coverings of a disk
$D^2$ with two branch points $A$ and $B$ in the interior of $D^2$.
These are determined by transitive representations $\rho$ of
$\pi_1 (D^2-\{A,B\})$ (which is free on the two meridian
generators, call them $x$ and $y$, pictured in Figure 13) into
$\Sigma_k$.

\begin{figure}
[ptbh]
\begin{center}
\includegraphics
{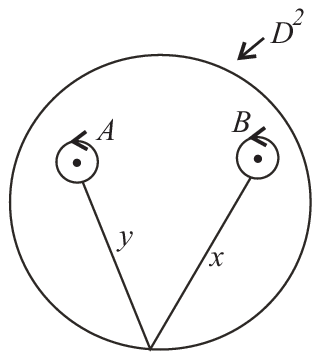}%
\caption{}%
\end{center}
\end{figure}

We are interested in a particular dihedral representation.  Let $P$
be a regular $k$--gon with its vertices labelled 1 through $k$ as in
Figure 14. We map $x$ to the reflection in the axis $l_1$ and $y$ to
the reflection in the axis $l_2$ of Figure 14.

\begin{figure}
[ptbh]
\begin{center}
\includegraphics
{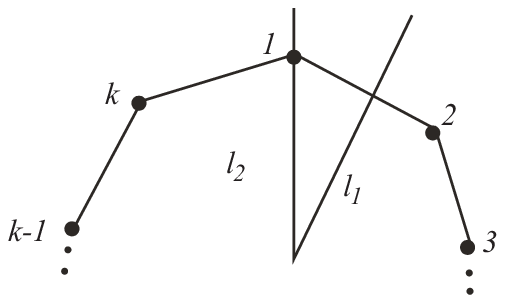}%
\caption{}%
\end{center}
\end{figure}

The reflections induce permutations of vertices and elements of
$\Sigma_k$.
 Then $\rho(x) = (1,2)(3,k)(4,k-1)\cdots$ and $\rho(y)=(2,k)(3,k-1)(4,k-2)\cdots$ are both products of disjoint transpositions and as $xy$ is a counterclockwise
rotation by $2\pi/k$, $\rho(xy)$ is the $k$--cycle $(1, k, k-1,
\cdots, 3,2)$.

Let $q: X \lra D^2$ be the branched covering induced by $\rho$.
Then $xy$ is a generator of $\pi_1(S^1)=Z$ where $S^1=\op{boundary }
D^2$ and restricting $q$ to boundary $X$ we see that $q: \partial X
\lra S^1$ is just the usual $k$--fold unbranched cover of $S^1$.
Using the branching data we compute the Euler characteristic of $X$
which is one.  Thus $\rho$ induces a $k$--fold branched cover
$q:D^{2} \lra D^2$ with branch set two points in interior $D^2$ and
branching type $\{1,2\}$.  The restriction of $q$ to $S^1 =\partial
X=\partial D^2$ is the usual $k$--fold unbranched covering of $S^1$
by $S^1$.  It is convenient to summarize the preceding in a
proposition.

\begin{proposition}\label{p2} Let $q_1: D^2 \lra D^2$ be the usual $k$--fold
cyclic covering of $D^2$; the one given by $q_1:z \lra z^n$ in
complex coordinates if $D^2$ is the unit disc in $C$.  Then there is
another $k$--fold branched covering $q:D^2\lra D^2$  such that the
branch set is two points in interior $D^2$, the branching is of type
$\{1,2\}$ and $q=q_1$ on $S^{1}=\op{boundary} D^2$.
\end{proposition}

By ``crossing" with $S^1$ the next proposition follows easily.

\begin{proposition}\label{p3} Let $q_1: S^1 \times D^2 \lra S^1 \times D^2$
be the usual $k$--fold cyclic covering of $S^1 \times D^2$ by $S^1
\times D^2$; the one given by $q_1:(e^{i\theta},z)\lra
(e^{i\theta}, z^n)$ in natural coordinates for $S^1 \times D^2$.
Then there is another $k$--fold branched covering $q:S^1 \times
D^2 \lra S^1 \times D^2$ such that the branch set equals $S^1
\times \{A,B\}$ where $A$ and $B$ are points in the interior of
$D^2$; the branching is of type $\{1,2\}$, and $q=q_1$ on boundary
$(S^1 \times D^2)$.
\end{proposition}

We return to our map $f_1$ which is an $m$--fold cyclic covering of
$S^3$ by $S^3$ with branch set and preimage of branched set the
$z$--axis.  We choose a natural solid torus neighborhood of the
$z$--axis in $S^3$ and its preimage and we coordinatize this
neighborhood so that $f_1$ is the usual $m$--fold cyclic covering of
$S^1 \times D^2$ by $S^1 \times D^2$ as in \propref{p3}.  This
neighborhood should be small enough so that it doesn't intersect the
link $L$ which is the branch set of the map $p:M^3\lra S^3$ defined
earlier. Then we use \propref{p3} to define a new map $f:S^3 \lra
S^3$ where outside the solid torus neighborhood $f=f_1$ and within
the solid torus neighborhood $f$ is like $q$ of \propref{p3}.  Thus
$f\circ p$ is a $3m$ to $1$ branched covering of $S^3$ by $M^3$ with
branch set $S^1 \times \{A,B\} \cup f (L)$.

The part of the branch set $f(L)=f_1(L)$ has one vertical component
and $n$--horizontal components and $n$--``peanut" components.  Via
an isotopy, if necessary, we may assume that the rotation $T_2$
leaves $S^1 \times \{A\}$ and $S^1 \times \{B\}$ invariant and that
its restriction to either component is just the usual $n$--fold
rotation.  The relevant part of the branch set for $f\circ p$ is
depicted in Figure 15. The branching is type $\{1,2\}$.

\begin{figure}
[ptbh]
\begin{center}
\includegraphics
{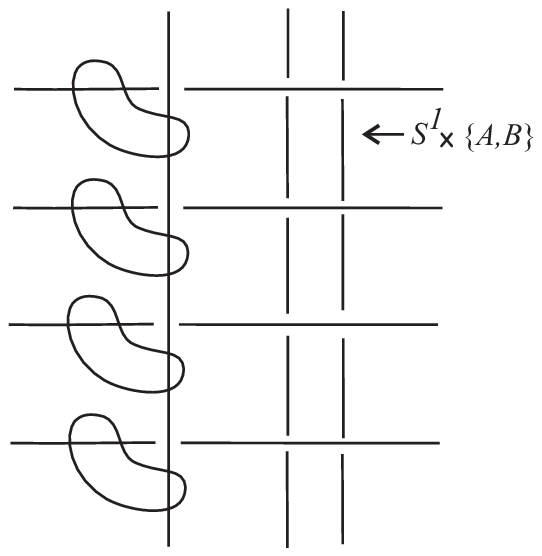}%
\caption{}%
\end{center}
\end{figure}

Next we let $g_1$ be the map $S^3\lra S^3/T_2 = S^3$.  Then $g_1
\circ f\circ p$ is a branched covering of $S^3$ by $M^3$ which is
$3mn$ to $1$ and has branch set equal to $g_1$(branch set ($f\circ
p$)) $\cup$ \{the circle $[z = 0,\, r = 2,\,\theta = \op{any}]$\}.
The branching is of type $\{1, 2\}$ on $g_1$(branch set $(f \circ
p)$) and of type $\{m\}$ on the circle $[z=0, r=z, \theta
=\op{any}]$.

We choose a solid torus neighborhood of the circle $z=0, r=z,
\theta=any$, small enough so that it does not intersect
$g_1$(branch set ($f\circ p$)), and so that it can be
coordinatized as $S^1 \times D^2$ with $g_1=q_1$ as in
\propref{p3}.

Then we define $g:S^3\lra S^3$ so that $g=g_1$ outside the torus
neighborhood and $g$ behaves like $q$ of \propref{p3} inside the
torus neighborhood.  The $3mn$ to $1$ branched covering $g\circ f
\circ p: M^3\lra S^3$ has branch set with three vertical components
and three horizontal components and one ``peanut" component as
pictured in Figure 16.  All branching is of type $\{1,2\}$.

\begin{figure}
[ptbh]
\begin{center}
\includegraphics
{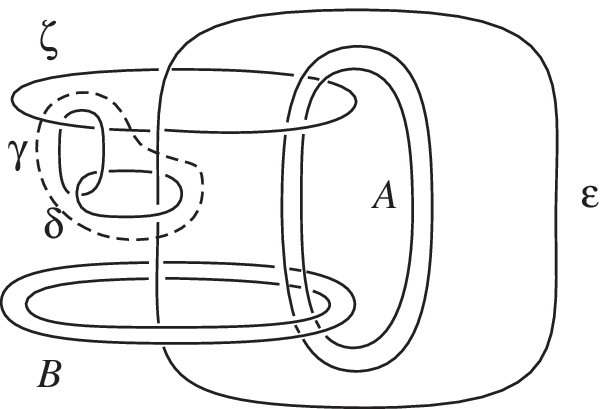}%
\caption{}%
\end{center}
\end{figure}

The link in Figure 16 can be isotoped to the link in Figure 17.  To
help the reader see this we have labelled the corresponding
components in Figures 16 and 17.  There are two obvious annuli in
Figures 16 and 17 labelled $A$ and $B$.  The other components are
labelled $\gamma,\delta,\epsilon,$ and $\zeta$.

\begin{figure}
[ptbh]
\begin{center}
\includegraphics
{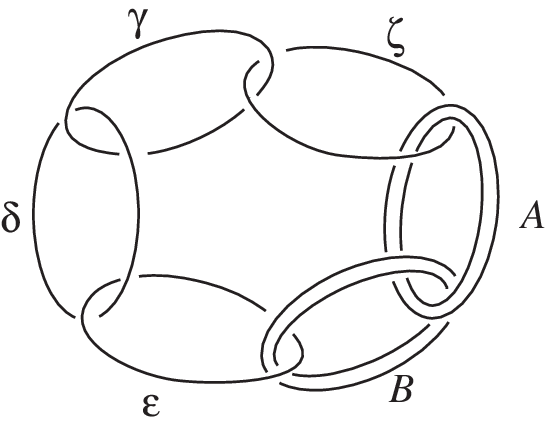}%
\caption{}%
\end{center}
\end{figure}

We may add four components, $\gamma_1,\delta_1,\epsilon_1, \zeta_1$,
to the link of Figure 17 so that there are annuli $C,D,E,$ and $F$
with boundaries $\gamma \cup \gamma_1$, $\delta \cup \delta_1$,
$\epsilon \cup \epsilon_1$ and $\zeta \cup \zeta_1$ respectively in
such a way that the new link has a $3$--fold rational symmetry. Let
$T_3$ be this $3$--fold rotation and let $h_1:S^3\lra S^3 = S^3/T_3$
be the resulting branched covering.  The map $h_1\circ g\circ f\circ
p$ is a $9mn$ to $1$ branched covering of $S^3$ by $M^3$ with
branching of type $\{1,2\}$ on the part of the branch set $h_1$
(branch set $g\circ f\circ p$) and branching of type $\{3\}$ on $h_1
(\op{axis} T_3)$.

This branch set is depicted in Figure 18.

\begin{figure}
[ptbh]
\begin{center}
\includegraphics
{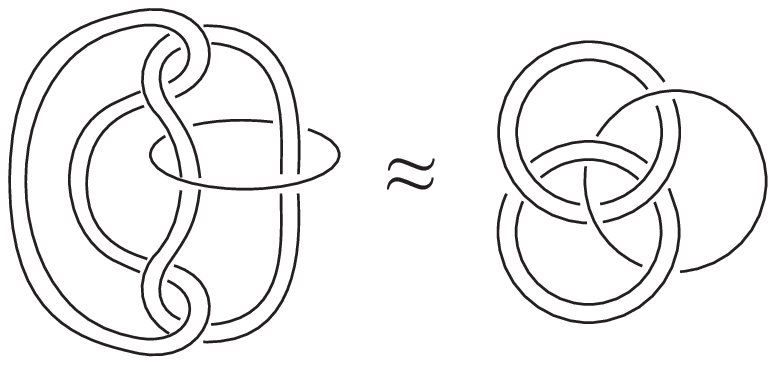}%
\caption{}%
\end{center}
\end{figure}

As before we replaced $f_1$ by $f$ and $g_1$ by $g$ we now replace
$h_1$ by $h$ using \propref{p3} where $h_1=h$ except in solid torus
neighbourhood of the axis of rotation of $h$ but the branching of
the map $h:S^3 \lra S^3$ is of type $\{1,2\}$.  The new branch set
is displayed in Figure 19.

\begin{figure}
[ptbh]
\begin{center}
\includegraphics
{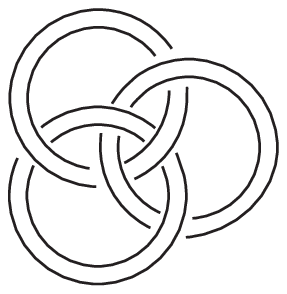}%
\caption{}%
\end{center}
\end{figure}

We call the branch set of Figure 19 the ``doubled Borromean rings".
We summarize this result in the form of a theorem \cite{LM1997}.

\begin{theorem}\label{t4} Let $M^3$ be a closed orientable $3$--manifold.
Then $M^3$ is a branched covering of $S^3$ with branch set the
doubled Borromean rings, and with branching of type $\{1,2\}$. That
is, the doubled Borromean rings are 2-universal.
\end{theorem}

If we use the right hand side of Figure 11 instead of the left hand
side we obtain the following theorem \cite{LM1997}.

\begin{theorem}\label{t5} Let $M^3$ be a closed orientable $3$--manifold.
Then $M^3$ is a branched covering of $S^3$ with branch set the
doubled Whitehead link and with all branching of type $\{1,2\}$.
That is, the doubled  Whitehead link is 2-universal.
\end{theorem}

The doubled Whitehead link is depicted in Figure 20.

\begin{figure}
[ptbh]
\begin{center}
\includegraphics
{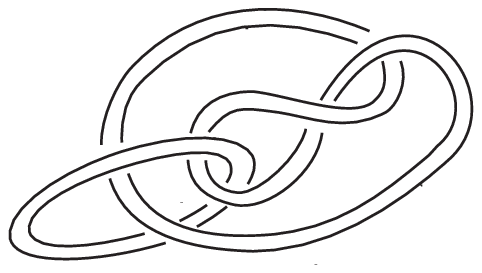}%
\caption{}%
\end{center}
\end{figure}

\thmref{t4} and \thmref{t5} show that the doubled Borromean rings
and doubled Whitehead link are 2-universal.  In (\cite{LM1997}) two
co--authors of this paper prove this result and additionally give
infinitely many examples of 2-universal links. The three component
link in the right bottom part of Figure 5.10 of \cite{LM1997} is a
minimal hyperbolic 2-universal link. In fact any proper sublink of
it is either a split link or a toroidal link.

We note, in passing, that 2-universal knots are known to exist
(\cite{HLM2004}) but so far there are no easy examples.

Our next task, which is the new idea of this paper, will be to
define a branched covering $k:S^3 \lra S^3$ with branch set the
Borromean rings for which the doubled Borromean rings occur as a
sublink of the preimage of the branch set.

We begin by tessellating $E^3$ by $2 \times 2 \times 2$ cubes all of
whose vertices have odd integer coordinates.  Let $\wh{U}$ be the
group generated by $180^{\circ}$ rotations in the axes $a,b,$ and
$c$ displayed in Figure 21.  The cube there is centered at the
origin.

\begin{figure}
[ptbh]
\begin{center}
\includegraphics
{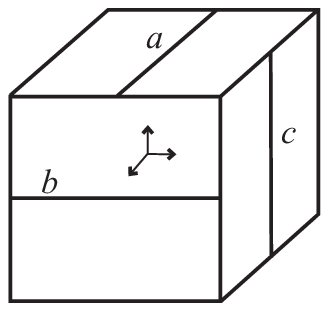}%
\caption{}%
\end{center}
\end{figure}

The group $\wh{U}$ is a well known Euclidean crystallographic group
that preserves the tessellation.  A fundamental domain for $\wh{U}$
is the $2 \times 2 \times 2$ cube of Figure 21 centered at the
origin.  The map $E^3 \lra E^3/\wh{U} \approx S^3$ is a branched
cover of $S^3$ by $E^3$ with branch set the Borromean rings.  This
gives $S^3$ the structure of a Euclidean orbifold with singular set
the Borromean rings and singular angle $180^{\circ}$.  We can see
that $E^3/\wh {U}$ equals $S^3$ with singular set the Borromean
rings by making face identifications in the fundamental domain of
Figure 21.  We do this in Figure 22.

\begin{figure}
[ptbh]
\begin{center}
\includegraphics
{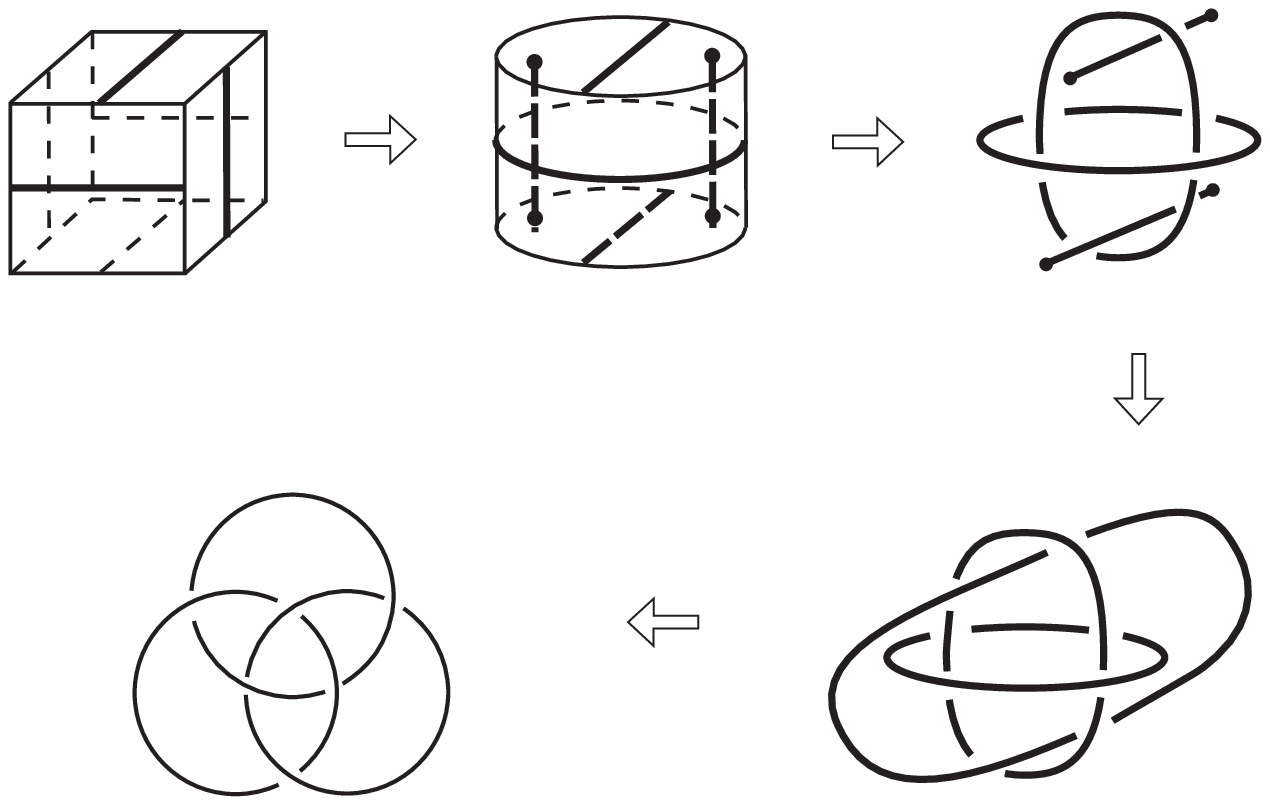}%
\caption{}%
\end{center}
\end{figure}

Next we consider a tessellation of $E^3$ by $6 \times 6 \times 6$
cubes with integer coordinates that are odd multiples of three.  Let
$\widetilde U$ be the group generated by $180^{\circ}$ rotations in
the axes $a'$, $b'$, $c'$ where $a' = (t,0,3)$, $ b'= (3,t,0)$ and
$c'= (0,3,t); -\infty <t<\infty$.  We can envision a fundamental
domain for $\widetilde{U}$ by looking at Figure 21 and imagining
prime accents  over $a,b,$ and $c$.  Of course $E^3/\widetilde U =
S^3$ and the map $E^3\lra E^3/\widetilde{U}$ is a branched covering
of $S^3$ by $E^3$.  As the rotations about $a'$, $b'$ and $c'$
belong to $\wh{U}$ we see that $\widetilde{U} \subset \wh{U}$ and we
see that $\big[\wh{U}:\widetilde{U}\big]=27$ by comparing the size
of fundamental domains. $\widetilde{U}$ is not a normal subgroup of
$\wh {U}$.  We are in fact really interested in the map $t:S^3 =
E^3/\widetilde{U} \lra E^3/\wh{U} = S^3$ induced by the inclusion of
$\widetilde{U}$ in $U$.

Consider the following commutative diagram

\[
\begin{CD}
E^3 @>id>> E^3 \\
@VVV  @VVV\\
S^3 = E^3/\widetilde U @>>> E^3/\widehat U = S^3\,.
\end{CD}
\]

The maps $E^3\lra E^3/\widetilde{U}$ and $E^3 \lra E^3/\wh{U}$ are
both branched covers of $S^3$ by $E^3$ with all branching of type
$\{2\}$.  We see, by considering images of meridian discs in
$E^3$, that the map $S^3 = E^3/\widetilde{U} \lra E^3/\wh{U} =
S^3$ is also a branched covering space map with branching of type
$\{1,2\}$.

The set of points in $S^3 = E^3/\widetilde{U}$ in the preimage of
the Borromean rings branch set in $E^3/\wh{U}=S^3$ for which the
branching of type $\{2\}$ is called the branch cover.  The set of
points in $S^3 = E^3/\widetilde{U}$ in the preimage of the Borromean
rings branch set in $E^3/\wh{U} = S^3$ for which the branching is of
type $\{1\}$ is called the pseudo branch cover.

We can compute the preimage of the branch set in $S^3 =
E^3/\widetilde{U}$ from a fundamental domain for $\widetilde{U}$,
which consists of $27$ $2 \times 2 \times 2$ cubes.

A point in this fundamental domain belongs to the branch cover if
and only if it belongs to an axis of rotation for $\wh{U}$ but
does not belong to an axis of rotation for $\widetilde{U}$.  A
point in this fundamental domain belongs to the pseudo branch
cover if and only if it belongs to an axis of rotation for
$\widetilde{U}$.  We do not need to compute the full preimage of
the Borromean rings in $S^3 = E^3/\widetilde{U}$, (which turns out
to be a $15$ component link), but only a certain sublink.

In Figure 23 we give a $6 \times 6 \times 6$ cube fundamental domain
for $\widetilde{U}$ and display only those axes of rotation of
$\wh{U}$ that lie in the faces of the cube.

The axes $a'$, $b'$ and $c'$ (and their analogues on the invisible
faces) are axes for $\widetilde{U}$ and so give rise to the pseudo
branch cover.  The axes $\widetilde a$, $\widetilde b$, and
$\widetilde c$ (and their analogues on the invisible faces) are axes
for $\wh{U}$, but not $\widetilde{U}$, and so give rise to a sublink
of the branch cover.

\begin{figure}
[ptbh]
\begin{center}
\includegraphics
{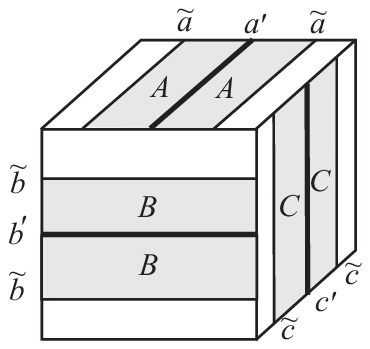}%
\caption{}%
\end{center}
\end{figure}

Also, axes $a'$ and $\widetilde {a}$ (resp. $b'$ and
$\widetilde{b}$; $c'$ and $\widetilde {c}$\,) lie on the boundary of
a rectangle $A$ (resp. $B$; $C$).  Something similar occurs on the
invisible faces. After identifications are made in the faces of the
cubes these rectangles become annuli.  In Figure 24 we make the
identifications and show that the doubled Borromean rings appear.

\begin{figure}
[ptbh]
\begin{center}
\includegraphics
{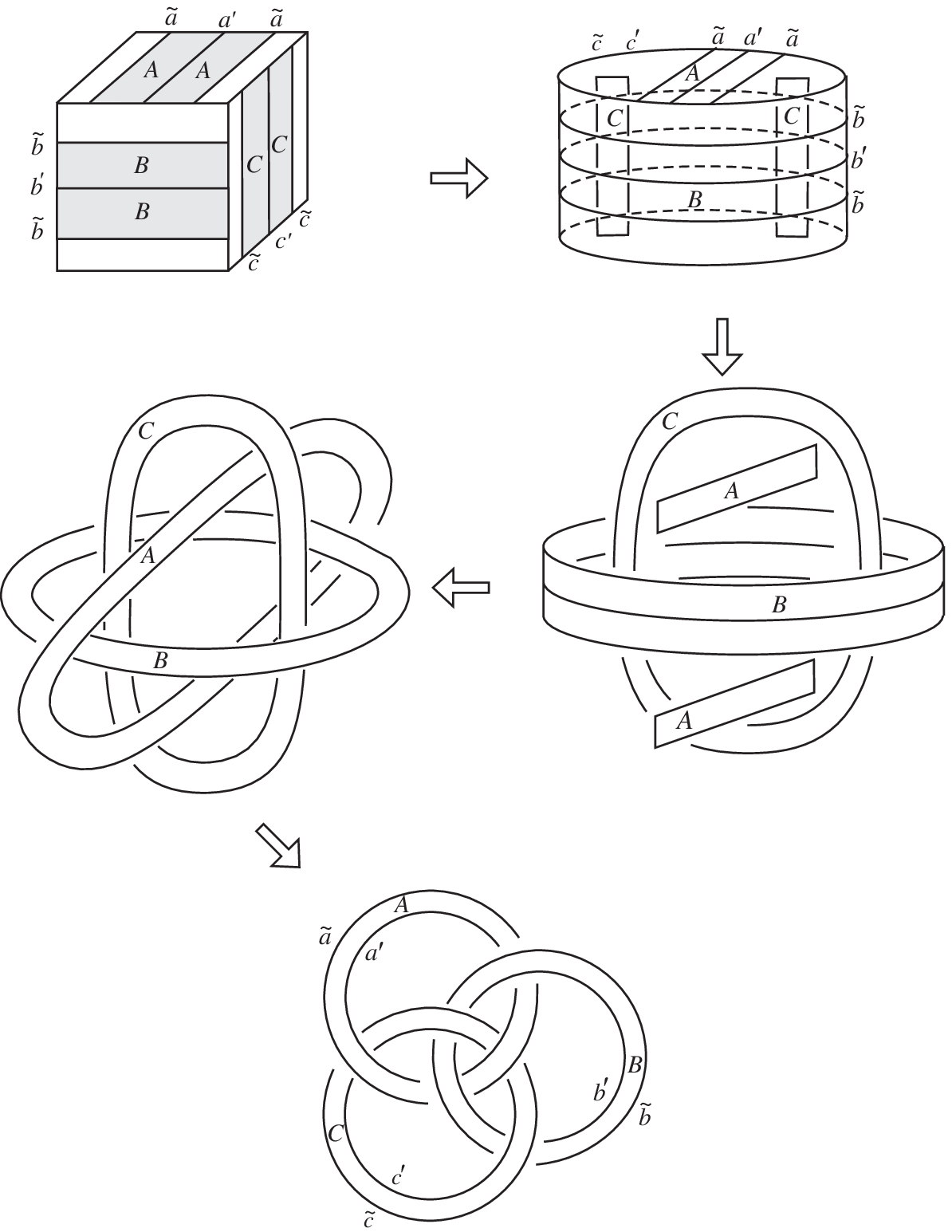}%
\caption{}%
\end{center}
\end{figure}

We summarize all this in the next proposition.

\begin{proposition}\label{p7} The map $t:S^3 = E^3/\widetilde{U} \lra
E^3/\wh{U} = S^3$ is a $27$ to $1$ irregular branched covering of
$S^3$ by $S^3$ with branch set the Borromean rings.

The doubled Borromean rings occur as a sublink of the preimage of
the branch set.

The doubled Borromean rings consist of three pairs of components.
Each pair bounds an annulus disjoint from the other pairs.  Each
pair is mapped to the same component of the Borromean rings.  And
each pair contains one component of the branch cover and one
component of the pseudo branch cover.
\end{proposition}
Thus far, starting with an arbitrary closed orientable $3$--manifold
we have defined a series of branched covering space maps.

\[
M^3\underset{p}{\overset {3-1}{\lra}} S^3 \overset{m-1}{\underset
{f}{\lra}} S^3 \overset{n-1}{ \underset{g}{\lra}} S^3 \underset{h}{
\overset{3-1}{\lra}} S^3 \underset{t}{ \overset{27-1}{\lra}}S^3\,.
\tag{1}
\]

The maps $f$ and $g$ only depend on $M^3$ in a superficial way; i.e.
they depend on the number of vertical and horizontal components in
the link that is the branch set for $p$, but not on the branching
itself and the maps $h$ and $t$ don't depend on $M^3$ at all.

In general, when one composes branch covering space maps
$a:X^3\lra Y^3$, $b:Y^3\lra Z^3$ one obtains a branch covering
space map $b\circ a: X^3 \lra Z^3$.

In the sequel we abbreviate branch set by $BS$.  Thus $BS(b \circ a)
= BS(b) \cup b(BS(a))$.  In our case $BS(f) \cap f(BS)(p) =
\emptyset$; $BS(g) \cap g(BS(f \circ p)) = \emptyset $; $BS(h) \cap
h(BS(g \circ f \circ p)) =\emptyset$.  And $h \circ g \circ f \circ
p$ is a branched covering space map of $S^3$ by $M^3$, of branching
type $\{1, 2\}$ with branch set the doubled Borromean rings.  Unlike
the previous three compositions $BS(t) = t(BS(h \circ g \circ f
\circ p))$ and the map $t \circ h \circ g \circ f \circ p$ is a
branched covering space map of $S^3$ by $M^3$ of branching type
$\{1, 2, 4\}$.

We summarize this as the next theorem

\begin{theorem}\label{t8} Let $M^3$ be a closed orientable $3$--manifold.
There is a branched covering space map $q:M^3\lra S^3$ with
branching type \{1, 2, 4\} and branch set the Borromean rings.

There is a $3$--fold simple branched coverings space map $p:M^3 \lra
S^3$ branched over a link, and there is a branched covering space
map $\phi:S^3\lra S^3$ with branch set the Borromean rings such that
$q=\phi\circ p$.
\end{theorem}
We remark that \thmref{t8} implies that the Borromean rings are
4-universal.  That is, that every closed orientable $3$--manifold
is a branched cover of $S^3$ with branch set the Borromean rings.
This was shown for the first time in \cite{HLM1983c}.  The proof
here refines the proof of \cite{HLM1983c} in the sense that bounds
are given on the branching.  The branching is of type \{1, 2, 4\}.
The refinement consists in introducing the map $t:S^3\lra S^3$
whose definition is based on the Euclidean orbifold structure of
$S^3$ with singular set the Borromean rings.

Part of \thmref{t8}, the first two sentences, was first proven in
(\cite{HLMW1987}).  But the branched covering space map $q:M^3\lra
S^3$ whose existence was proven there did not factor as in the rest
of the statement in \thmref{t8}.  This factorization is crucial to
the coming proof of \thmref{t9}, the principal result of this paper.

The sphere $S^3$ has a hyperbolic orbifold structure with singular
set the Borromean rings and singular angle $90^{\circ}$.  The maps
$t$, $h$, $g$, $f$, and $p$ are used to pull back the hyperbolic
orbifold structure on $S^3$ to hyperbolic orbifold structures on
$S^3,S^3,S^3,S^3$, and $M^3$ respectively.  Thus there is a sequence
of groups, orbifold groups, $G \subset G_1 \subset G_2 \subset G_3
\subset G_4\subset U$ such that $M^3 = H^3/G$, $S^3 = H^3/G_1$, $S^3
= H^3/G_2$, $S^3 = H^3/G_3$, $S^3 = H^3/G_4$, $S^3 = H^3/U$ and the
following diagram is commutative.  The vertical arrows are
isomorphisms.
\[
\begin{CD}
H^3/G @>>> H^3/G_1 @>>> H^3/G_2 @>>> H^3/G_3 @>>> H^3/G_4 @>>> H^3/U \\
@VVV   @VVV   @VVV   @VVV   @VVV   @VVV \\
M^3 @>>p> S^3 @>>f> S^3 @>>g> S^3 @>>h> S^3 @>>t> S3
\end{CD}
\]

In particular we see that $[G_1:G]=3$.

If $T$ is any rotation contained in $G_1$ but not contained in $G$
then $G_1 = \langle G, T\rangle$, as $[G_1:G]=3$ and 3 is prime.
 We summarize these
remarks in the main theorem of this paper.

\begin{theorem}\label{t9} {(Geometric branched covering space theorem.)} Let
$M^3$ be a closed orientable $3$--manifold.  Then there are
subgroups $G$ and $G_1$ of the universal group $U$ such that
$[G_1:G]=3$ and $[U:G]<\infty$ and $M^3 = H^3/G$ and $S^3 =
H^3/G_1$.

The map induced by the inclusion of groups $H^3/G \lra H^3/G_1$ is a
$3$--fold simple branched covering of $S^3$ by $M^3$.
\end{theorem}

We recall (\cite{HLMW1987}) that if $M^3$ is as in the statement of
\thmref{t9} then $\pi_1(M^3) \cong G/\op{TOR}(G)$ where
$\op{TOR}(G)$ is the subgroup of $G$ generated by rotations. In
particular $M^3$ is simply connected if and only if $G=\op{TOR}(G)$.
Then $G_1$ is generated by $G$ and any one rotation not in $G$.

An interesting problem is to classify the finite index subgroups of
$U$ that are generated by rotations.  We intend to return to this
theme in a future paper.

Applying the theory of  associated regular coverings to the above
situation we obtain an interesting property on the involved groups
that  restricts their study to a subclass of the class of subgroups
of the universal group $U$ defining the same variety. Next we
explain this.

 In general, given a
covering $p: M \lra N$ branched over $L$, with monodromy $\omega
:\pi_{1}(N\setminus L)\lra \Sigma_{n}$, the \emph{ associated
regular covering} $q: X \lra N$ is the branched covering determined
by the monodromy $\eta \circ \omega :\pi_{1}(N\setminus L)\lra
\Sigma_{\sharp Im(\omega )}$, where $\eta$ is the regular
representation of the group $Im(\omega )$. Recall that $q = u \circ
p$ where $u: X \lra M$ is a regular (branched or unbranched)
covering. Actually $(q|_{X \setminus q^{-1}(L)})_{\star}(\pi_{1}(X
\setminus q^{-1}(L))=Ker(\omega )$.

The monodromy of $p: M=H^{3} / G \lra S^{3}=H^{3} / G_{1} $ is a
homomorphism $\omega :\pi_{1}(S^{3} \setminus L)\lra \Sigma_{3}$
where the image of every meridian element of $L$ is a transposition
$(i,j)$, $1\leq i <j \leq 3$. Therefore, $Im(\omega )=\Sigma_6$, and
the monodromy $\eta \circ \omega :\pi_{1}(S^{3} \setminus L)\lra
\Sigma_{6}$ sends every meridian element of $L$ to the product of
three different transpositions. Thus,  $u:X\lra M$ is a 2-fold
covering branched over the pseudobranch cover of $p: M=H^{3} / G
\lra S^{3}=H^{3} / G_{1} $. The covering $q = u \circ p :X\lra
S^{3}=H^{3} / G_{1} $ is a regular 6-fold covering branched over L
with all branching indexes equal to 2. Actually, the map $u$ can be
used to pull back the hyperbolic orbifold structure on $M$, so that
there exist a normal subgroup $G_{0}\triangleleft G$, such that
$[G:G_{0}]=2$ and \[u:X=H^{3}/ G_{0} \lra M=H^{3} / G\] is an
orbifold covering. Observe that $G_{0}$ is a normal subgroup of
$G_{1}$.

The following diagram of orbifold coverings is commutative, where
$G_{0}\subset G'\subset G_{1}$ and $[G_{1}:G']=2$, $[G':G_{0}]=3$.

\[
\begin{CD}
H^3/G_{0}  @>u>{2:1}> M=H^3/G  \\
  @V{3:1}V{p'}V  @V{p}V{ 3:1}V  \\
  H^3/G' @>{u'}>{2:1}> S^{3}=H^3/G_{1}
\end{CD}
\]

The covering $p'$ is unbranched and $u'$ is the cyclic 2-fold
covering branched over $L$. Observe that $G_{0}$ is a normal
subgroup of $G_{1}$. We summarize these remarks in the following
theorem.

\begin{theorem}\label{t10} Let $M^3$ be a closed orientable $3$--manifold.
Let $G$ and $G_1$ be the subgroups of the universal group $U$ given
in \thmref{t9}. Then there exist a subgroup $G_{0}$ of index 2 of
$G$ such that $G_{0}$ is a normal subgroup of $G_{1}$.
\end{theorem}

Recall that every finite index subgroup $G$ of the universal group
$U$ gives rise to a 3-manifold $M=H^{3}/G$, but infinitely many
finite index $G$'s produce the same manifold. The  above theorem
restricts the class of subgroups of $U$ to consider in order to
construct all  closed 3-manifolds.

\end{section}

\enddocument